\documentclass[11pt]{article}
\usepackage{mathrsfs}
\usepackage{latexsym,lineno}
\usepackage{epsfig}
\usepackage{color}
\usepackage{amsmath}\usepackage{fleqn}\usepackage{verbatim}\usepackage{epsf}
\usepackage{amsthm}\usepackage{graphicx, float}\usepackage{graphicx}
\usepackage{amsfonts}\usepackage{amssymb}\usepackage{graphpap}
\usepackage{epic}\usepackage{curves}
\newtheorem{thm}{Theorem}

\newtheorem{obs}{Observation}

\newtheorem{lem}{Lemma}

\title{On the super connectivity of Kronecker products of graphs\thanks{Research was partially supported by
PuJiang Project of Shanghai (No. 09PJ1405000) and Shanghai Leading
Academic Discipline Project (No. S30104).}}
\author{Hechao Wang,  \,  Erfang Shan\\
 {\small Department of
Mathematics, Shanghai University, Shanghai 200444, P.R. China} }
\date{}
\begin{document}

\maketitle
\begin{abstract}
Let $G_1$ and $G_2$ be two  graphs. The {\em Kronecker product}
$G_1\times G_2$ has vertex set $V(G_1\times G_2)=V(G_1)\times
V(G_2)$ and  edge set $E(G_1\times
G_2)=\{(u_1,v_1)(u_2,v_2):u_1u_2\in E(G_1)$ and $v_1v_2\in
E(G_2)\}$. A graph $G$ is {\em super connected}, or simply {\em
super}-$\kappa$, if every minimum separating set is the neighbors of
a vertex of $G$, that is, every separating set isolates a vertex. In
this paper we show that for  an arbitrary graph $G$ with
$\kappa(G)=\delta(G)$ and $K_n$ $(n\geq3)$ a complete graph on $n$
vertices, $G\times K_n$ is super-$\kappa$, where $\kappa(G)$ and
$\delta(G)$ are the connectivity and the minimum degree of $G$,
respectively.
\bigskip



\end{abstract}

%=================================================

\section{Introduction and terminology}
%We consider only finite and simple graphs (without loops or multiple
%edges).
%For notation and terminology not defined here we refer to
%Bondy \cite{Bondy}. Specifically,
Throughout this paper, only undirected simple graphs without loops
or multiple edges are considered. Unless otherwise stated, we follow
Bondy \cite{Bondy} for terminology and definitions.

Let $G=(V(G),E(G))$ be a graph.
 For two vertices $u,v\in V(G)$, $u$ and $v$ are {\em
neighbors} if $u$ and $v$ is adjacent.
 The set of vertices adjacent
to the vertex $v$ is called the {\em neighborhood} of $v$ and
denoted by $N(v)$, i.e., $N(v) =\{u\mid uv\in E(G)\}$. The {\em
degree} of $v$ is equal to  $|N(v)|$, denoted by $d_G(v)$ or simply
$d(v)$. The number $\delta(G)= {\rm min} \{ d_G(v) \mid v\in V(G)
\}$ is the {\em minimum degree} of $G$. For a subset $S\subseteq
V(G)$, the neighborhood of $S$ is $N(S)=\bigcup_{v\in S}N(v)$. The
subgraph induced by $S$ is denoted by $G[S]$, and let $d_{S}(v)$
denote the number of vertices in $S$ that are adjacent to the vertex
$v$. As usual, $K_n$ denotes the complete graph on $n$ vertices and
$C_n$ is the cycle on $n$ vertices.

A set $S\subset V$ is a {\em separating set} of a connected graph
$G$, if either $G-S$ disconnected or reduces to the trivial graph
$K_1$. The {\em connectivity} of  $G$, denoted by $\kappa(G)$, is
the minimum cardinality of a separating set of $G$. In particular,
 $\kappa(K_n)=n-1$ and
 $\kappa(G)=0$ if and only if $G$ is  disconnected or a $K_1$.
Clearly,   $\kappa(G)\le \delta(G)$. A graph $G$ with minimum degree
$\delta(G)$ is {\em maximally connected} if $\kappa(G) = \delta(G)$.

The notion of super--connectedness proposed in
\cite{Boesch1,Boesch2,Boesch3} aims at pushing the analysis of the
connectivity properties of graphs beyond the standard connectivity.
A graph $G$ is {\em super connected}, or simply {\em
super-$\kappa$}, if every minimum separating set is the neighbors of
a vertex of $G$, that is, every separating set isolates a vertex.
Observe that a super-connected graph $G$ is necessarily maximally
connected, i.e., $\kappa(G)=\delta(G)$, but the converse is not
true. It is easily to see from the cycle graph $C_n\ (n\geq6)$.

%An {\em homomorphism} between two graphs $G$ and $H$ is a mapping
%$\varphi$ from $V(G)$ to $V(H)$ such that $\{u,v\}$ is an edge of
%$G$, then $\{\varphi(u),\varphi(v)\}$ is an edge of $H$. We say that
%$\varphi$ is an {\em isomorphic} if $\varphi$ is bijective and
%$\varphi^{-1}$ is also an homomorphism.

The Kronecker product, together with the Cartesian, the strong, and
the lexicographic product, is one of the four standard graph
products \cite{ik}.
%The Cartesian product $G_1\Box G_2$ is defined
%as: $V(G_1\Box G_2)=V(G_1)\times V(G_2)$ and $E(G_1\Box
%G_2)=\{(u_1,v_1)(u_2,v_2):(v_1=v_2$ and $u_1u_2\in E(G_1))$ or
%$(u_1=u_2$ and $ v_1v_2\in E(G_2)\}$.
The Kronecker product of two graphs $G_1$ and $G_2$ is defined as
$V(G_1\times G_2)=V(G_1)\times V(G_2)$ and $E(G_1\times
G_2)=\{(u_1,v_1)(u_2,v_2):u_1u_2\in E(G_1)$ and $v_1v_2\in E(G_2)\}$
(see, \cite{Bottreou, Weichsel}). The Kronecker product has been
introduced and studied from several points of view and is known
under many different names, for instance as the direct product,
cardinal product, categorical product, tensor product and cross
product. Moreover, it is universal in the sense that every graph is
an induced subgraph of a suitable direct product of complete graphs
\cite{ne}. The Kronecker product of graphs has been extensively
investigated concerning graph colorings, graph recognition and
decomposition, graph embedding, matching theory, stability and
domination theory in graphs (see, for example, \cite{Alon, Bresal,
me}). The properties on the structure of Kronecker product of graphs
can be found in \cite{Guji, super-k, Mamut, spacapan2}. One has
known that it has many interesting applications, for instance it can
be used in automata theory \cite{gh}, complex networks \cite{lck}
and modeling concurrency in multiprocessor systems \cite{lb}.

Miller \cite{mi} and Weichsel \cite{Weichsel} investigated the
connectedness of Kronecker product of two connected graphs.
Recently, the connectivity of Cartesian products and strong products
of two connected graphs have been studied, and in all cases the
explicit formulae have been obtained in terms of the
 graph invariants of the factor graphs (see,
\cite{spacapan1, spacapan2,Xu} for more details). The connectivity
of Kronecker products of graphs seems to be more complex than that
with the Cartesian or strong products.
%In \cite{Mamut} Mamut and Vumar gave the connectivity of Kronecker product of two complete
%graphs.
Guji and Vumar \cite{Guji} presented the connectivity of Kronecker
product of a bipartite graph and a complete graph and they proposed
to investigate the connectivity of Kronecker product of a nontrivial
graph and a complete graph. Recently, Wang and Xue \cite{Wang}
settled the problem and they obtained the following result.
\begin{thm}[\cite{Wang}]\label{yibantu}
For any graph $G$, $\kappa(G\times K_n)= \min \{n\kappa(G),
(n-1)\delta(G)\}$ for $n\geq3$.
\end{thm}

%\begin{cor}[\cite{Wang}]\label{coro}
%If $G$ is a tripartite graph with $\kappa(G)=\delta(G)$, then
%$\kappa(G\times K_n)=\delta(G\times K_n)=(n-1)\delta(G)$ for
%$n\geq3$.
%\end{cor}

Recently,  Guo et al. \cite{super-k} studied the super connectivity
of Kronecker product of a bipartite graph and a complete graph and
they proved the following result.

\begin{thm}[\cite{super-k}]\label{thm2}
If $G$ is a bipartite graph with $\kappa(G)=\delta(G)$, then
$G\times K_n\ (n\ge3)$ is super-$\kappa$.
\end{thm}

%\begin{thm} [\cite{Guji}]\label{thm1}
%If $G$ is a bipartite graph, then $\kappa(G\times K_n)=$ min
%$\{n\kappa(G), (n-1)\delta(G)\}$ for $n\geq3$.
%\end{thm}
%In \cite{Guji} the authors also proposed the following conjecture.
%\begin{conjecture}\label{conj1}
%If $G$ is a nontrivial graph, then $\kappa(G\times K_n)=$ min
%$\{n\kappa(G), (n-1)\delta(G)\}$ for $n\geq3$.
%\end{conjecture}
In this paper, motivated by the above results, we further
investigate the super connectivity of Kronecker product of an
arbitrary  graph and a complete graph $K_n$ ($n\ge 3$). Our main
result is as follows.
\begin{thm}\label{main}
For an arbitrary graph $G$ with $\kappa(G)=\delta(G)$, $G\times K_n\
(n\ge3)$ is super-$\kappa$.
\end{thm}

%\begin{corollary}\label{coro}
%If $G$ is a tripartite graph with $\kappa(G)=\delta(G)$, then
%$\kappa(G\times K_n)=\delta(G\times K_n)=(n-1)\delta(G)$ for
%$n\geq3$.
%\end{corollary}
\section{Preliminaries}

In this section we give some properties on Kronecker product of
graphs, and we will use them in the proof of our main result.
%From the definition, the following observation due to \cite{Bottreou} is immediate.
We first give two known results.
\begin{obs} [\cite{Bottreou}] \label{lem1}
Let $H=G_1\times G_2=(V(H), E(H))$. Then

$(1)$  $|V(H)|=|V(G_1)|\cdot |V(G_2)|$,

$(2)$  $|E(H)|=2|E(G_1)|\cdot |E(G_2)|$,

$(3)$  for every $(u,v)\in V(H)$, $d_H((u,v))=d_{G_1}(u)\cdot
d_{G_2}(v)$.
\end{obs}
By Observation \ref{lem1},   we have $\delta(G\times
K_n)=(n-1)\delta(G)$ for any graph $G$.

\begin{lem}[\cite{Weichsel}] \label{connected} Let $G_1$ and $G_2$ be
connected graph. The graph $H=G_1\times G_2$ is connected if and
only if $G_1$ or $G_2$ contains an odd cycle.
\end{lem}

When considering the Kronecker product of a graph $G$ and $K_n$
$(n\geq3)$, we shall always let $V(G)=\{u_1,u_2,\ldots,u_m\}$,
$V(K_n)=\{v_1,v_2,\ldots,v_n\}$ and set $S_i=\{u_i\}\times V(K_n)$,
for $i=1,2,\ldots,m$. Then $S_i$ is an independent set in $G\times
K_n$, and $V(G\times K_n)$ has a partition $V(G\times K_n)=S_1\cup
S_2\cup\ldots\cup S_m$.

Let $S\subseteq V(G\times K_n)$ satisfy the following three
conditions:

$(1).\ |S|=(n-1)\delta(G)$, and

$(2).\ S_i':=S_i-S\neq \varnothing$, for $i=1,2,\ldots,m$, and

$(3).\ G\times K_n-S$ has no isolated vertex.

 Associated with $G,\ S$ and $S_i'$, we
define the following new graph $G^*$ as described in \cite{Wang}.

${\rm (i)}$. $V(G^*)=\{S_1',S_2',\ldots, S_m'\}$, and

${\rm (ii)}$. $E(G^*)=\{S_i'S_j':\ E(S_i',S_j')\neq \varnothing\}$,
where $E(S_i',S_j')$ denotes the collections of all edges in
$G\times K_n-S$ with one end in $S_i'$ and the other in $S_j'$.

Next we give two lemmas about the connectedness of $G^*$ and the
structure of $S_i'$, which play a key role in the proof of our main
result. In the next proofs, we always assume
$\kappa(G)=\delta(G)>0$.

 In
\cite{Wang} proved that $G^*$ is connected if $G$ is connected and
$|S|<(n-1)\delta(G)$. In fact, when $|S|=(n-1)\delta(G)$ here,  we
still get the  result by using  the same method in \cite{Wang}.
\begin{lem}\label{GSTar}
If $G$ is connected, then $G^*$ is connected.
\end{lem}

\begin{lem}\label{Si}
Let $G$ be nonbipatite graph with $\kappa(G)=\delta(G)$. Then for
any vertex of $G^*$, $S_i'$, as a subset of $V(G\times K_n)$, it is
contained in the vertex set of some component of $G\times K_n-S$.
\end{lem}

\noindent{\bf Proof.} It suffices to prove the lemma for $i=1$. We
consider the cardinality of the set $S_i'$.

If $|S_1'|=1$, then the assertion holds trivially. So we may assume
that $|S'|\ge 2$. We  consider the following two cases:

{\bf Case 1:} $|S_1'|\geq3$. Suppose to the contrary that $S_1'$ is
not contained in any component of $G\times K_n-S$. Then there must
exist a component, say $C$, such that $0< |S_1'\cap V(C)|\leq
|S_1'|/2<|S_1'|-1$. Let $(u_1,v_p)\in S_1'\cap V(C)$. By the
conditions $(1)$ and $(3)$ of the definition of $S$,
 $G\times K_n-S$ has at least one vertex, say $(u_j,v_q)(j\in \{2,\ldots,m\})$,
such that $(u_j,v_q)\ {\rm and}\ (u_1,v_p)$ are neighbors. Clearly,
$(u_j,v_q)\in V(C)$, and $S_1'-\{(u_1,v_q)\}\subseteq V(C)$ since
every vertex in $S_1'-\{(u_1,v_q)\}$ is adjacent to $(u_j,v_q)$. It
follows $|S_1'\cap V(C)|\geq |S_1'|-1$, a contradiction.

{\bf Case 2:} $|S_1'|=2$. Let $Z^*=\{S_j':j\in\{1,2,\ldots,m\},\
|S_j'|=1\}$ and $C^*$ be a component of $G^*-Z^*$ which containing
$S_1'$. Let $|C^*|=r$, without loss of generality, we may assume
$V(C^*)=\{S_1',S_2',\ldots,S_r'\}$.

Since each $S_j'\in V(C^*)$ contains at least two elements, any edge
$S_j'S_k'$ in $C^*$ implies every vertex in $S_k'$ has at least one
neighbor in $S_j'$ in $G\times K_n-S$. Therefore, if there is a
vertex $S_j'$ in $C^*$ contains in the vertex set of some component
$C$ of $G\times K_n-S$, then every $S_k'$ with $S_k'S_j'\in E(C^*)$
is contained in $V(C)$ as each $|S'_j|>1$. It follows from the
connectedness of $C^*$ that $\bigcup_{i=1}^rS_i'\subseteq V(C)$ and
hence $S_1'\subseteq V(C)$.

By Case 1, we may assume each $S_j'\in V(C^*)$ contains exactly two
elements. Let $S_j'=\{u_j\}\times F_j$, $j=1,2,\ldots,r$, and
$F_j\subseteq V(K_n)$.

Suppose that there exists an edge $S_j'S_k'$ in $C^*$ with $F_j\neq
F_k$. It is easily to see that $S_j'\cup S_k'$ induces a connected
subgraph of $G\times K_n$. This implies that $S_j'$ and $S_k'$ are
contained in the same component, say  $C$, of $G\times K_n-S$. As
mentioned above, we have $S_1'\subseteq\bigcup_{i=1}^rS_i'\subseteq
V(C)$,  the lemma follows.

Otherwise,  by the connectedness of $C^*$, we have
$F_1=F_2=\cdots=F_r$. Notice that $C^*$ and $G[\bigcup_{i=1}^rS_i']$
are isomorphic to $G[u_1,u_2,\ldots,u_r]$ and
$G[u_1,u_2,\ldots,u_r]\times K_2$, respectively. If we can prove
that $G[u_1,u_2,\ldots,u_r]\times K_2$ is connected in $G\times
K_n-S$, then the lemma follows.

Note that $G$ is connected. If $G[u_1,u_2,\ldots,u_r]$ contains an
odd cycle, then the assertion follows by Lemma \ref{connected}.

Suppose that $G[u_1,u_2,\ldots,u_r]$ does not contain an odd cycle,
that is, $G[u_1,u_2,\ldots,u_r]$ is bipartite.
This implies that $\delta(G[u_1,u_2,\ldots,u_r])\leq
r/2$. Let $j\in\{1,2,\ldots,r\}$ such that
$d_{G[u_1,u_2,\ldots,u_r]}(u_j)\\=\delta(G[u_1,u_2,\ldots,u_r])$.

Let $u_k$ be a neighbor of $u_j$ in $G$. Then either $S_k'\in Z^*$
or $S_k'$ is  adjacent to $S_j'$ in $C^*$. Thus,
\begin{equation}\label{1}
\delta(G)\leq d_G(u_j)\leq
d_{C^*}(S_j')+|Z^*|=d_{G[u_1,u_2,\ldots,u_r]}(u_j)+|Z^*|\leq
r/2+|Z^*|.
\end{equation}
Therefore, we have
\begin{equation}\label{2}
(n-1)\delta(G)=|S|\geq (n-2)r+(n-1)|Z^*|\geq
(n-1)\Big(\frac{r}{2}+|Z^*|\Big)\geq (n-1)\delta(G).
\end{equation}
This means that the equations holds in (1) and (2). Hence, $n=3$,
$\delta(G)=d_G(u_j)=r/2+|Z^*|$ and
$d_{G[u_1,u_2,\ldots,u_r]}(u_j)=\delta(G[u_1,u_2,\ldots,u_r])=r/2$.
Since $G[u_1,u_2,\ldots,u_r]$ is bipartite, it is $r/2$-regular, so
each vertex $u_j$ ($1\le j\le r$) has the same degree $r/2+|Z^*|$ in
$G$. This implies that each $S_j'$ ($1\le j\le r$) is adjacent to
all the vertices of $Z^*$ in $G^*$.

We claim that $Z^*\neq \varnothing$.

Indeed, if $Z^*=\varnothing$,
then $C^*=G^*$, so $G[u_1,u_2,\ldots,u_r]=G[u_1,u_2,\ldots,u_m]=G$,
which contradicts our assumption that $G$ is a nonbipartite graph.
Clearly, $G[u_1,u_2,\ldots,u_r]\times K_2$ can be connected by the
vertices of $Z^*$, as desired.  $\Box$

By the definition, the following lemma is straight.
\begin{lem}\label{4}
Let $m=|G|\leq2$ and $u_i$ be any vertex of $G$. Then

{\rm (1)}. $\delta(G-u_i)\geq \delta(G)-1$, and

{\rm (2)}. $\kappa(G-u_i)\geq \kappa(G)-1$.
\end{lem}

\section{Proof of the main result}
Now we are ready to give the proof of Theorem \ref{main}. By Theorem
\ref{2}, we only need to show that the assertion in Theorem
\ref{main} is true for a nonbipartite graph $G$. Therefore, we
always assume that $G$ is nonbipartite in our proof below.

%\noindent{\bf Remark.} Our proof follows the lines of Theorem 2 in
%\cite{Wang}, with modifications needed because of the super
%connectivity.

\noindent{\bf Proof of Theorem \ref{main}.}  If $G$ is disconnected,
i.e., $\kappa(G)=\delta(G)=0$, then $G\times K_n$ is disconnected,
and the assertion holds. So we may assume that $G$ is connected and
$\delta(G)\geq1$.

If we can show that for every subset $S$ of $G\times K_n$ with
$|S|=(n-1)\delta(G)$, either $G\times K_n-S$ is connected or
$G\times K_n-S$ has an isolated vertex, then we are done. Our proof
is by contradiction. Suppose that $G\times K_n$ is not
super-$\kappa$. Then there is a separating set $S$ with
$|S|=(n-1)\delta(G)$ such that $G\times K_n-S$ is not connected but
has no isolated vertex. If we can show that $G\times K_n-S$ is
connected, then we shall arrive at a contradiction and the assertion
holds.

We will distinguish two possibilities as follows.

{\bf Case 1}: If $S$ satisfies condition (2), i.e., $S_i':=S_i-S\neq
\varnothing,$ for $i=1,2,\ldots,m$. It follows from Lemma
\ref{GSTar} and Lemma \ref{Si} that $G\times K_n-S$ is connected.

{\bf Case 2}: $S$ does not satisfies condition (2). Then there
exists an $S_i$ contained in $S$. So $S-S_i\subseteq V((G-u_i)\times
K_n)$ and
\begin{eqnarray*}
|S-S_i|&=&|S|-n=(n-1)\delta(G)-n\\
&<& \min \{n\kappa(G)-n,(n-1)\delta(G)-(n-1)\}\\
&\leq&\min \{n\kappa(G-u_i),(n-1)\delta(G-u_i)\},
\end{eqnarray*}
the last inequality above follows from Lemma \ref{4}.

By Theorem \ref{yibantu}, we have $\kappa((G-u_i)\times K_n)=\min
\{n\kappa(G-u_i),(n-1)\delta(G-u_i)\}$, Hence, $(G-u_i)\times
K_n-(S-S_i)$ is connected. Note that $(G-u_i)\times
K_n-(S-S_i)=G\times K_n-S$. Hence $G\times K_n-S$ is connected.

In all cases, we show that $G\times K_n-S$ is connected, this
contradicts our assumption on $S$. So the assertion follows.

%If $\kappa(G)=\delta(G)$, then by Theorem \ref{main}, we have
%$\kappa(G\times K_n)=\delta(G\times K_n)$. By making use of the same
%method, we can give the super connectivity of Kronecker products of
%a tripartite graph and a complete graph.

%In the proof of the main theorem, we use one proper classify to the
%elements of the  matrix $(a_{ij})_{3\times 3}$. we can also use the
%method in this paper to prove next theorem only need to analysis
%some extreme conditions.

%\begin{thm}\label{thm super-k}
%If $G$ is a tripartite graph with $\kappa(G)=\delta(G)$, then
%$G\times K_n$ $(n\geq3)$ is super-$\kappa$.
%\end{thm}

%But we can't use the same method to get the similar category in
%$(a_{ij})_{4\times4}$. If we can get a proper classify use other
%method, we can prove the next conjecture use the idea in this paper.

%\begin{conjecture}\label{conjecture2}
%If $G$ is a $r$-partite graph, then $\kappa(G\times K_n)=$ min
%$\{n\kappa(G), (n-1)\delta(G)\}$ for $n\geq r$.
%\end{conjecture}

%\bigskip
%{\bf \Large Acknowledgements}

%{\rm The authors are grateful to the referees for their valuable
%comments, which have led to improvements in the presentation of the
%paper.}

\end{document}